\documentclass[12pt]{amsart}
\usepackage[final]{epsfig}
\usepackage{graphics}
\usepackage{float}
\usepackage{amsmath}
\usepackage{amsfonts}
\usepackage{latexsym}
\usepackage{amssymb}
\usepackage{graphicx}
\usepackage{url}
\usepackage{epstopdf}
\usepackage{verbatim}
\usepackage[margin=1in]{geometry}
\usepackage{amsthm}
\usepackage{enumitem}
\usepackage{tikz}
\usepackage[stable]{footmisc}

\makeatletter
\newtheorem*{rep@theorem}{\rep@title}
\newcommand{\newreptheorem}[2]{%
\newenvironment{rep#1}[1]{%
 \def\rep@title{#2 \ref{##1}}%
 \begin{rep@theorem}}%
 {\end{rep@theorem}}}
\makeatother

\newtheorem{lemma}{Lemma}[section]
\newtheorem{proposition}[lemma]{Proposition}

\newtheorem{theorem}[lemma]{Theorem}
\newtheorem{definition}[lemma]{Definition}
\newtheorem{corollary}[lemma]{Corollary}

\newtheorem*{theorem*}{Theorem}
\newreptheorem{theorem}{Theorem}
\newtheorem{question}{Question}[section]

\newcommand{\proofend}{$\Box$\bigskip}

\newcommand{\N}{{\mathbb N}}
\newcommand{\R}{{\mathbb R}}

\def\proof{\paragraph{Proof.}}

\DeclareMathOperator{\Gram}{Gram}
\DeclareMathOperator{\GramS}{Gram_{\succeq 0}}

\newcommand{\CommentE}[1]{{\color{blue} \sf ($\diamondsuit$ #1 $\diamondsuit$)}}

\begin{document}

\title{Dimension of Gram Spectrahedra of Univariate Polynomials}

\author{Emmanuel Tsukerman
}

\date{\today}

\address{Department of Mathematics, University of California,
Berkeley, CA 94720-3840}
\email{e.tsukerman@berkeley.edu}

\maketitle
\setcounter{tocdepth}{1}

\begin{abstract}
The Gram Spectrahedron of a polynomial parametrizes its sums-of-squares representations. In this note, we determine the dimension of Gram Spectrahedra of univariate polynomials. 
\end{abstract}

  \section{Notation}

  \noindent
  $P_{n,2d} \text{ is the cone of nonnegative polynomials of degree no more than }2d \text{ in } n \text{ variables.}$
  
  \noindent
  $P_{n,2d}^{>0} \text{ is the cone of positive polynomials of degree no more than }2d \text{ in } n \text{ variables.}$
 
 \noindent
  $\Sigma_{n,2d} \text{ is the cone of sum-of-squares polynomials of degree no more than }2d \text{ in } n \text{ variables.}$ 
 
\noindent
$S_N(\R) \text{ is the vector space of } N \times N \text{ real symmetric matrices}.$   
 
 \noindent
$P_N^+(\R) \text{ is the cone of } N \times N \text{ positive semidefinite real symmetric matrices}.$

\section{Introduction}

The problem of understanding the relationship between nonnegative polynomials and sums-of-squares of polynomials goes back to Hilbert and his predecessors. Clearly a sum-of-squares of polynomials is nonnegative. However, the converse holds in general only in the following special cases:

\begin{theorem} (Hilbert) \cite[Pg. 59]{MR3050241} \label{Hilbert} The cones $\Sigma_{n,2d}$ and $P_{n,2d}$ are equal if and only if 
\begin{enumerate}
\item $n=1$ (univariate polynomials).
\item $2d=2$ (quadratic polynomials).
\item $n=2, 2d=4$ (bivariate quartics).
\end{enumerate}
\end{theorem}

Sums-of-squares questions for polynomials can fruitfully be  approached via semidefinite optimization theory. Indeed, let $f$ be a (multivariate) polynomial of degree $2d$ and let $X$ be the vector of monomials of degree up to $d$. It is well-known that
\[
f \text{ is a sum of squares} \iff \exists \text{ } Q \succeq 0: f=X^t Q X.
\]

\begin{definition}
\normalfont The \textit{Gram matrices of $f$} is the set 
\[
\Gram(f):=\{Q \in S_N(\R):f=X^t Q X\}
\]
The spectrahedron
\[
\GramS(f):=\{Q \in \Gram(f): Q \succeq 0\}
\]
is called the \textit{Gram spectrahedron} of $f$.
\end{definition}

In particular, $f$ is a sum-of-squares if and only if its Gram Spectrahedron is nonempty.

 Recent work includes \cite{MR2781949}, in which the Gram Spectrahedron of a real ternary quartic is studied, and \cite{QuarticSpectrahedra}, in which it is shown that the Gram Spectrahedron of a univariate sextic is affinely isomorphic to a spectrahedron associated to a Kummer surface.

\section{Main Results}

The following theorem gives a complete description of the dimension of the Gram Spectrahedron of a univariate polynomial.

\begin{theorem}\label{univariateDimension}
Let $f \in P_{1,2d}$ and write
\[
f(x)=(x-r_1)^{2e_1}(x-r_2)^{2e_2}\cdots (x-r_k)^{2e_k} g(x)
\]
for $r_i \in \R, e_i \in \N$ and $g(x) \in P_{1,2(d-e)}^{> 0}$, where $e=\sum_i e_i$. Then
\[
\GramS(f) \simeq \GramS(g).
\]
In particular,
\[
\dim \GramS(f)=\dim \GramS(g)=\binom{d-e}{2}.
\]
\end{theorem}

\section{Generalities}

\begin{lemma} \label{FullDimensionalFiber}
Let $f: \R^n \rightarrow \R^d$ be a surjective linear map,
$K$ a full-dimensional convex set in $\R^n$ and 
$L = f(K)$ its image in $\R^d$. Then every point $q$ that lies in the interior of $L$ has a preimage $p$ that lies in the interior of $K$.
\end{lemma}

\proof

 Let $A$ be the matrix for $f$. By the Singular Value Decomposition, we may view $A$ as $UDV$ with $U,V$ orthogonal. Thus we can reduce to the situation when $f$ is just a projection: $f(x_1,\ldots,x_n)=(x_1,\ldots,x_d)$. In addition, by considering $f$ to be the composition of a series of maps of the same type, $\R^n \xrightarrow{f_{1}} \R^{n-1} \xrightarrow{f_{2}} \R^{n-2} \xrightarrow{f_{3}} \ldots \xrightarrow{f_{n-d}} \R^{d}$, we reduce to $d=n-1$.

Take $C_\delta(q)=\{x \in \R^{n-1}: \|x-q\|^2 < \delta^2\} \subset L$.  The preimage of the ball under $A$ is 
\[
C_\delta(q):=\{x \in \R^n: \|Ax-q\|^2 < \delta^2\},
\]
a full-dimensional cylinder in $\R^n$.
Since $B_\delta(q) \subset L$, for each $z$ satisfying $\|z-q\|^2 < \delta^2$ there exists some $g(z_1,\ldots,z_{n-1})$ such that
\[
(z_1,\ldots,z_{n-1},g(z_1,\ldots,z_{n-1})) \in K.
\]

By convexity, we may take $g$ to be a continuous function. Let 
\[
G=\{(z_1,\ldots,z_{n-1},g(z_1,\ldots,z_{n-1}))\}.
\] 
Since $K$ is full-dimensional, it must contain a ball $B$. Considering the convex hull of $B$ and $G$ shows that there is a preimage of $q$ in the interior of $K$.
\proofend

 \begin{proposition}
  Let $f \in \Sigma_{n,2d}^o$ be a polynomial of degree $2d$ in $n$ variables lying in the interior of the sum of squares cone. Then the Gram Spectrahedron of $f$ is full-dimensional:
  \[
  \dim \GramS(f)=\binom{\binom{d+n}{n}+1}{2}-\binom{2d+n}{n}.
  \]
  \end{proposition}
  
  \proof
  Let $N:=\binom{d+n}{n}$.
  There is a surjective linear map $\phi: S_N(\R) \rightarrow \R_{\leq 2d}[x_1,x_2,\ldots,x_n]$ given by
\[
\phi(Q)=X Q X^t,
\]
where $X$ is the vector consisting of the monomials in $x_1,x_2,\ldots,x_n$ of degree up to $d$. The cones $\Sigma_{n,2d}$ and $P_N^+(\R)$ are proper \cite[Theorem 3.26]{MR3050241}.
Since $f$ lies in the interior of the cone $\Sigma_{n,2d}$, by Lemma \ref{FullDimensionalFiber}, there is a point in the fiber $\phi^{-1}(f)$ which lies in the interior of $P_N^+(\R)$. The dimension of $S_N(\R)$ is $\binom{N+1}{2}$ and the dimension of $\R_{\leq 2d}[x_1,\ldots,x_n]$ is $\binom{n+2d}{2d}$. The result follows.
  \proofend

\begin{corollary} \label{fullDim}
If $f$ is in 
\begin{enumerate}[label=(\roman*)]
\item $P_{1,2d}^{>0}$ 
\item $P_{n,2}^{>0}$
\item $P_{2,4}^{>0}$,
\end{enumerate}
 then the Gram Spectrahedron of $f$ is full-dimensional:
\begin{enumerate}[label=(\roman*)]
\item $\dim \GramS(f)= \binom{d}{2}$
\item $\dim \GramS(f)= 0$
\item $\dim \GramS(f)=6$.
\end{enumerate}
\end{corollary}

\proof
In these cases, we have $\Sigma_{n,2d}=P_{n,2d}$ by Theorem \ref{Hilbert}. The boundary of $P_{n,2d}$ is determined by the vanishing of the discriminant. Since $f$ has no real roots, it lies in the interior of $P_{n,2d}=\Sigma_{n,2d}$.
\proofend

\section{Dimension of Gram Spectrahedra of Univariate Polynomials}

In this section, we determine the dimension of the Gram Spectrahedron of a univariate polynomial.

\begin{lemma} \label{gramIso}
Let $f(x) \in \R_{2d}[x]$ and let $a,b \in \R$. Let $A:=aR+bS \in \R^{(d+2) \times (d+1)}$ with
\[
R=\left( \begin{array}{cccccccc}
0 & 0 & \cdots & 0 \\
1 & 0 & \cdots & 0 \\
0 & 1 & \ddots  & \vdots\\
\vdots & \ddots & \ddots & 0\\
0 &  \cdots & 0 & 1  
   \end{array} \right)
\]
and
\[   
    S=\left( \begin{array}{cccccccc}

1 & 0 & \cdots & 0 \\
0 & 1 & \ddots  & \vdots\\
\vdots & \ddots & \ddots & 0\\
0 &  \cdots & 0 & 1  \\
0 & 0 & \cdots & 0 
   \end{array} \right)
   \]
Then the map
\[
\phi:\GramS(f(x)) \rightarrow \GramS((ax+b)^2f(x))
\]
\[
M \rightarrow AMA^t
\]
induced by multiplication by $(ax+b)^2$ is an affine isomorphism.
\end{lemma}

\proof
If $\GramS(f(x))$ is empty then clearly so is $\GramS((ax+b)^2f(x))$ and vice versa. So assume otherwise. 
Let 
\[
f(x)= \sum_{i} q_i(x)^2
\]
be a sum of squares representation. Multiplying by $(ax+b)^2$ gives a sum of squares representation of $(ax+b)^2f(x)$:
\[
(ax+b)^2f(x)=\sum_i ((ax+b)q_i(x))^2.
\]
Conversely, given a sum of squares representation of $(ax+b)^2f(x)$,
\[
(ax+b)^2f(x)=\sum_i q_i(x)^2,
\]
we evaluate both sides at the root $r$ of $ax+b$. This shows that $q_i(r)=0$ for each $r$. Dividing both sides by $(ax+b)^2$ yields a sum of squares representation for $f(x)$:
\[
f(x)=\sum_i (\frac{q_i(x)}{ax+b})^2.
\]
These operations are clearly inverses of one another. 

 At the level of the Gram Spectrahedron, we have the following.  Set $(c_0\text{ }c_1\text{ }c_2\text{ }\dots\text{ }c_d)^t \in \R^{d+1}$ to be the vector of coefficients of $q_i$.
The vector $q'_i$ of coefficients of $(ax+b)q_i^t x$ is then equal to $aRq+bq=(aR+bS)q$. Then $q' q'^t=(aR+bS)q q^t (aR+bS)$. Therefore the image of $Q$ is $AQA^t$. This map is clearly linear and preserves the property of being positive semidefinite.
\proofend

\proof (of Theorem \ref{univariateDimension})
The result follows from Corollary \ref{fullDim} and Lemma \ref{gramIso}.

\proofend

The dimension of Gram Spectrahedra of quadratics is trivially zero. 

\begin{question}
Find the dimension of Gram Spectrahedra of other families of polynomials, such as bivariate quartics.
\end{question}

\bigskip
{\bf Acknowledgments}.  
This material is based upon work supported by the National Science Foundation Graduate Research Fellowship under Grant No. DGE 1106400. Any opinion, findings, and conclusions or recommendations expressed in this material are those of the authors(s) and do not necessarily reflect the views of the National Science Foundation.

\bibliographystyle{alpha}
\bibliography{bibliography}

\end{document}